\documentstyle[amstex,amssymb,11pt, amscd, epsfig]{article}

\def\e0{\epsilon_{0}}

\newtheorem{theorem}{Theorem}[section]
\newtheorem{prop}[theorem]{Proposition}
\newtheorem{lemma}[theorem]{Lemma}

\begin{document}

\date{\empty}
\title{{\bf Counting Elliptic Curves in K3 Surfaces \vskip.2in}}
\author{By Junho Lee and Naichung Conan Leung \\
School of Mathematics, University of Minnesota}
\maketitle

\begin{abstract}
We compute the genus $g=1$ family GW-invariants of K3 surfaces for
non-primitive classes. These calculations verify G\"{o}ttsche-Yau-Zaslow
formula for non-primitive classes with index two. Our approach is to use the
genus two topological recursion formula and the symplectic sum formula to
establish relationships among various generating functions.
\end{abstract}


\addtocounter{section}{0}


\vskip.15in 

\vskip.4in

The number of elliptic curves in K3 surfaces $X$ representing a homology
class $A\in H_{2}\left( X,{\Bbb Z}\right) $ and pass through one generic
point depends only on the self-intersection number $A\cdot A=2d-2$ and the
index\footnote{%
The index of $A$ is the largest positive integer $r$ such that $r^{-1}A$ is
integral. An index one class is called primitive.} $r$ of the class $A$ (
\cite{bl1}, \cite{bl2}). We denote it as $N_{1}\left( d,r\right) $. The
conjectural formula of G\"{o}ttsche \cite{go} for elliptic curves, which
generalize the Yau-Zaslow formula \cite{yz}, assert that the generating
function for those numbers $N_{1}(d,r)$'s for any given index $r$ is given
by
\begin{equation}
\sum_{d\geq 0}N_{1}(d,r)\,t^{d}\,=\,t\,G_{2}^{\prime }(t)\,\prod_{l\geq 1}%
\Big(\frac{1}{1-t^{l}}\Big)^{24}  \label{GYZ-for}
\end{equation}
where $G_{2}(t)$ is the Eisenstein series of weight 2, i.e.
\[
G_{2}(t)\,=\, \sum_{d\geq 0}\sigma(d)\,t^{d} \ \ \ \mbox{where}\ \ \
\sigma(d)\,=\,\sum_{k|d}k,\ d\geq 1\ \ \ \mbox{and}\ \ \ \sigma(0)\,=\,-%
\tfrac{1}{24}.
\]
In particular $N_{1}\left( d,r\right) $ should be independent of the index
of the homology class. In \cite{bl1,l2}, this formula was verified for
primitive classes by using modified Gromov-Witten invariants of K3 surfaces.
In this article, we verify the formula (\ref{GYZ-for}) for index two classes
by computing the family GW-invariants defined in \cite{l1}. Our main theorem
is the following result.

\smallskip

\begin{theorem}
\label{T:Main} Let $X$ be a K3 surface and $A/2\in H_{2}\big(X;{\Bbb Z}\big)$
be a primitive class. Then, the genus $g=1$ family GW-invariant of $X$ for
the class $A$ is given by
\begin{equation}
GW_{A,1}^{{\cal H}}\,=\,GW_{B,1}^{{\cal H}}\,+\,2\,GW_{A/2,1}^{{\cal H}}
\label{0.2}
\end{equation}
where $B$ is any primitive class with $B^{2}=A^{2}$.
\end{theorem}

\bigskip

To explain the equivalence between the above theorem and the
G\"{o}ttsche-Yau-Zaslow formula, we first notice that the {\em family
Gromov-Witten }invariant counts the number of $J$-{\em holomorphic maps} for
any complex structure $J$ in the twistor family of the K3 surface $X$. As
explained in \cite{bl1}, there is a unique $J$ in the twistor family which
supports holomorphic curves representing $A$ and this justifies the use of
our family invariant. An important issue is the distinction between
holomorphic maps to $X$ and holomorphic curves in $X$. This is because a
multiple curve in $X$ can be the image of different holomorphic maps. This
issue does not arise when the homology class $A/2$ they represent is
primitive and therefore $N_{1}\left( d^{\prime },1\right) =GW_{A/2,1}^{{\cal %
H}}$, where $\left( A/2\right) ^{2}=2d^{\prime }-2$. The number of genus one
holomorphic maps covering a fix elliptic curve with degree $r$ equals the
partition function $p\left( r\right) $, for instance $p\left( 2\right)
=1+2=3 $. Each primitive elliptic curve in $X$ contributes $3$ to $GW_{A,1}^{%
{\cal H}}$, but it only contributes $1$ to $N_{1}\left( d,2\right) $.
Therefore the above theorem is equivalent to,
\[
N_{1}(d,2)=N_{1}(d,1)\text{.}\,
\]
Together with the validation of the formula for primitive classes, this
implies the G\"{o}ttsche-Yau-Zaslow formula for the number of elliptic
curves in K3 surfaces representing index two homology classes.

\bigskip

The organization of this paper is as follows: The construction of family
GW-invariants is briefly described in section 1. This section also contain a
family version of the composition law. In section 2, using the composition
law and the genus $g=2$ TRR (Topological Recursion Relation) formula \cite{g}
we establish the $g=2$ TRR formula for family GW-invariants. In section 3,
we prove Theorem~\ref{T:Main} by combining that TRR formula with the
symplectic sum formulas of \cite{ll}.

\bigskip \noindent {\bf Acknowledgments\thinspace :} The first author would
like to thank Thomas Parker for his extremely helpful discussions and he is
also grateful to Eleny Ionel, Bumsig Kim and Ionut Ciocan-Fontaine for their
useful comments. In addition, the first author wish to thank Ronald
Fintushel for his interest in this work and especially for his
encouragement. The second author is partially supported by NSF/DMS-0103355.

\vskip 1cm


\section{Composition Law for Family GW-Invariants}

\label{section1} \bigskip

This section briefly describes family GW-invariants defined in \cite{l1}.
Let $X$ be a K\"{a}hler surface with a K\"{a}hler structure $(\omega,J,g)$.
For each 2-form $\alpha$ in the linear space
\[
{\cal H}=\mbox{Re}\big(H^{2,0}\oplus H^{0,2}\big)
\]
we define an endormorphism $K_{\alpha}$ of $TX$ by the equation $\langle
u,K_{\alpha}v\rangle=\alpha(u,v)$. Since $Id+JK_{\alpha}$ is invertible,
\[
J_{\alpha} = (Id + JK_{\alpha})^{-1}J\,(Id + JK_{\alpha})
\]
is an almost complex structure on X.

Denote by $\overline{{\cal F}}_{g,k}(X,A)$ the space of all stable maps $%
f:(C,j)\rightarrow X$ of genus $g$ with $k$-marked points which represent
the homology class $A$. For each such map, collapsing unstable components of
the domain determines a point in the Deligne-Mumford space $\overline{{\cal M%
}}_{g,k}$ and evaluation of marked points determines a point in $X^{k}$.
Thus we have a map
\begin{equation}
\overline{{\cal F}}_{g,k}(X,A)\,@>\ st\times ev\ >>\,\overline{{\cal M}}%
_{g,k}\times X^{k}  \label{st-ev}
\end{equation}
where $st$ and $ev$ denote the stabilization map and the evaluation map,
respectively. On the other hand, there is a {\em generalized orbifold bundle}
$E$ over $\overline{{\cal F}}_{g,k}(X,A)\times {\cal H}$ whose fiber over $%
\big((f,j),\alpha \big)$ is $\Omega _{jJ_{\alpha }}^{0,1}(f^{\ast }TX)$.
This bundle has a section $\Phi $ defined by
\[
\Phi \big(f,j,J_{\alpha }\big)\,=\,df+J_{\alpha }\,df\,j\,.
\]
When $X$ is a K3 surface and $A\neq 0$, the moduli space $\Phi ^{-1}(0)=%
\overline{{\cal M}}_{g,k}^{{\cal H}}\big(X,A\big)$ is compact. By the same
manner as in the theory of the ordinary GW-invariants \cite{lt}, this
section then gives rise to a well-defined rational homology class
\[
GW_{g,k}^{{\cal H}}(X,A)\,\in \,H_{2r}\big(\overline{{\cal F}}_{g,k}(X,A);%
{\Bbb Q}\big)\ \ \ \mbox{where}\ \ r=g+k.
\]
The family GW invariants of $(X,J)$ are defined by
\[
GW_{g,k}^{{\cal H}}(X,A)\big(\beta ;\alpha \big)\,=\,GW_{g,k}^{{\cal H}%
}(X,A)\,\cap \,\big(st^{\ast }(\beta ^{\ast })\,\cup \,ev^{\ast }(\alpha
^{\ast })\big)
\]
where $\beta ^{\ast }$ and $\alpha ^{\ast }$ are Poincar\'{e} dual of $\beta
\in H_{\ast }\big(\overline{{\cal M}}_{g,k};{\Bbb Q}\big)$ and $\alpha \in
H_{\ast }\big(X^{k};{\Bbb Q}\big)$, respectively.

\smallskip In \cite{l1} the first author proved that the above family
GW-invariants of K3 surfaces are same as the invariants defined by Bryan and
Leung \cite{bl1} using the twistor family. In particular, they are
independent of complex structures and for any two classes $A,B$ of the same
index with $A^{2}=B^{2}$, we have
\begin{equation}
GW_{g,k}^{{\cal H}}(X,A)\,=\,GW_{g,k}^{{\cal H}}(X,B).  \label{T:BL}
\end{equation}
\ \medskip We will often denote the above family GW-invariants simply as $%
GW_{A,g}^{{\cal H}}$.

\medskip

The family GW-invariants have a property analogous to the composition law of
ordinary GW-invariants (cf. \cite{rt}). Consider a node of a stable curve $C$
in the Deligne-Mumford space $\overline{{\cal M}}_{g,k}$. When the node is
separating, the normalization of $C$ has two components. The genus and the
number of marked points decompose as $g=g_{1}+g_{2}$ and $k=k_{1}+k_{2}$ and
there is a natural map
\[
\sigma : \overline{{\cal M}}_{g_{1},k_{1}+1}\times\overline{{\cal M}}%
_{g_{2},k_{2}+1}\to \overline{{\cal M}}_{g,k}.
\]
defined by identifying $(k_{1}+1)$-th marked points of the first component
to the first marked point of the second component. We denote by $PD(\sigma)$
the Poincar\'{e} dual of the image of this map $\sigma$. For non-separating
node, there is another natural map
\[
\theta:\overline{{\cal M}}_{g-1,k+2}\to \overline{{\cal M}}_{g,k}
\]
defined by identifying the last two marked points. We also write $PD(\theta)$
for the Poincar\'{e} dual of the image of this map $\theta$.

\smallskip Recall that the ordinary GW-invariants of $K3$ surfaces are all
zero except for trivial homology class.

\begin{prop}[\protect\cite{l1}]
\label{cl} Let $\{H_{a}\}$ be a base of $H_{\ast }(X;{\Bbb Z})$ and $%
\{H^{a}\}$ be its dual base with respect to the intersection form.

\begin{enumerate}
\item[{\rm (a)}]  Given any decomposition $g=g_{1}+g_{2}$ and $k=k_{1}+k_{2}$%
, we have
\begin{align*}
& GW_{A,g}^{{\cal H}}(PD(\sigma );\alpha _{1},\cdots ,\alpha _{k}) \\
& =\ \sum_{a}GW_{A,g_{1}}^{{\cal H}}(\alpha _{1},\cdots ,\alpha
_{k_{1}},H_{a})\,GW_{0,g_{2}}(H^{a},\alpha _{k_{1}+1},\cdots ,\alpha _{k}) \\
& +\ \sum_{a}GW_{0,g_{1}}(\alpha _{1},\cdots ,\alpha
_{k_{1}},H_{a})\,GW_{A,g_{2}}^{{\cal H}}(H^{a},\alpha _{k_{1}+1},\cdots
,\alpha _{k})
\end{align*}
where $GW_{0,g_{1}}$\ and $GW_{0,g_{2}}$ denotes the ordinary GW invariants
of $K3$ surfaces.

\item[(b)]  $GW_{A,g}^{{\cal H}}(PD(\theta );\alpha _{1},\cdots ,\alpha
_{k})=\sum_{a}GW_{A,g-1}^{{\cal H}}(\alpha _{1},\cdots ,\alpha
_{k},H_{a},H^{a}).$
\end{enumerate}
\end{prop}

\vskip 1cm


\section{Topological Recursion Relations}

\bigskip

The first composition law, Proposition~\ref{cl}\thinspace a, relates family
invariants and ordinary invariants (for trivial homology class) of K3
surfaces. In this section, we first recall the ordinary GW-invariants of a
closed symplectic 4-manifold for the trivial homology class. As in \cite{l2}%
, we then combine the composition law with the genus $g=2$ TRR formula \cite
{g} to establish a family version of $g=2$ TRR formula.

\bigskip

Let $\tau_{i}$ be the first Chern class of the line bundle ${\cal L}_{i}\to
\overline{{\cal M}}_{g,k}^{{\cal H}}\big(X,A\big)$ whose geometric fiber at
the point $\big(C;x_{1},\cdots,x_{k},f,\alpha\big)$ is $T^{*}_{x_{i}}C$.

\begin{lemma}
\label{GW=0} Let $X$ be a closed symplectic 4-manifold. Its Gromov-Witten
invariants satisfy the following properties:

\begin{enumerate}
\item[(a)]  $GW_{0,0}\big(\alpha _{1},\cdots ,\alpha _{k}\big)=0$ unless $%
k=3 $ and $\sum \,{\rm deg}\alpha _{i}=4$\thinspace {\rm . In that case,} $%
GW_{0,0}\big(\alpha _{1},\alpha _{2},\alpha _{3}\big)=\int_{X}\alpha
_{1}^{\ast }\,\alpha _{2}^{\ast }\,\alpha _{3}^{\ast }$,

\item[(b)]  $GW_{0,1}\big(\alpha _{1},\cdots ,\alpha _{k}\big)=0$ unless $%
k=1 $\thinspace {\rm . In that case} $GW_{0,1}\big(\alpha \big)=-\frac{1}{24}%
\,c_{1}(\alpha )$,

\item[(c)]  $GW_{0,2}\big(\tau (\alpha _{1}),\alpha _{2}\big)=0$. When $%
c_{1}\left( T_{X}\right) =0$, we also have $GW_{0,2}\big(\tau (\alpha )\big)%
=0$.
\end{enumerate}
\end{lemma}

\noindent {\bf Proof. } (a) and (b) directly follows from Proposition 1.4.1
of \cite{km1}. On the other hand, the formula (7) of \cite{km1} says that
\[
GW_{0,2}\big(\tau (\alpha _{1}),\alpha _{2}\big)\,=\,\big(\alpha _{1}\cdot
\alpha _{2}\big)\int_{\overline{{\cal M}}_{2,2}}\lambda _{2}^{2}\,\psi
_{1}\,,
\]
and
\[
GW_{0,2}\big(\tau (\alpha )\big)\,=\,-c_{1}(\alpha )\int_{\overline{{\cal M}}%
_{2,1}}\lambda _{1}\,\lambda _{2}\,\psi _{1}
\]
where $\lambda _{i}=c_{i}(E)$ is the Chern class of the Hodge bundle $E$ and
$\psi _{i}=c_{1}(L_{i})$ is the first Chern class of the line bundle $%
L_{i}\rightarrow \overline{{\cal M}}_{g,k}$ whose geometric fiber at the
point $\big(C;x_{1},\cdots ,x_{k}\big)$ is $T_{x_{i}}^{\ast }C$. The first
invariant in (c) is zero since $\lambda _{2}^{2}=0$ (cf. \cite{m}), while
the second one vanishes when $c_{1}\left( T_{X}\right) =0$. \nopagebreak%
\hskip.1in {\ $\Box $ }\penalty10000 \hskip\parfillskip
\vskip0.5cm

Let $E(2)\rightarrow {\Bbb P}^{1}$ be an elliptic K3 surface with a section
of self intersection number $-2$. Denote by $s$ and $f$ the section class
and the fiber class, respectively. We will also denote by $(S,F)$ either $%
(2s,f)$ or $(s-3f,2f)$.

\smallskip

\begin{prop}
\label{P:prop} The family Gromov-Witten invariants of an elliptic K3 surface
$E\left( 2\right) $ satisfy the following formula,
\[
GW_{S+dF,2}^{{\cal H}}\big(\tau _{1}(F),\tau _{2}(F)\big)=\,\ -\frac{2}{3}%
\,GW_{S+dF,1}^{{\cal H}}\big(pt\big)\,+\,\frac{(d-2)^{2}\,}{9}\,GW_{S+dF,0}^{%
{\cal H}}\,.
\]
\end{prop}

\noindent {\bf Proof. } In the same manner as for ordinary GW-invariants,
combining Proposition~\ref{cl} with the formula (5) of \cite{g} yields an
expression for the family invariants. Let $\{H_{a}\}$ and $\{H^{a}\}$ be
bases of $H^{*}\big(E(2);{\Bbb Z}\big)$ which are dual by the intersection
form. We then have
\begin{align*}
&\ \, GW^{{\cal H}}_{S+dF,2}\big(\tau_{1}(F),\tau_{2}(F)\big)  \nonumber \\
\ =&\, \sum_{a}2\,\Big( GW^{{\cal H}}_{S+dF,2}\big(\tau_{1}(F),H_{a}\big)\,
GW_{0,0}\big(H^{a},F\big)\, +\, GW_{0,2}\big(\tau_{1}(F),H_{a}\big)\, GW^{%
{\cal H}}_{S+dF,0}\big(H^{a},F\big) \Big)  \nonumber \\
\ -&\, \sum_{a,b} GW^{{\cal H}}_{S+dF,0}\big(F,H_{a}\big)\, GW_{0,0}\big(%
F,H_{b}\big)\, GW_{0,2}\big(H^{a},H^{b}\big)  \nonumber \\
-&\, \sum_{a,b} GW_{0,0}\big(F,H_{a}\big)\, GW^{{\cal H}}_{S+dF,0}\big(%
F,H_{b}\big)\, GW_{0,2}\big(H^{a},H^{b}\big)  \nonumber \\
-&\, \sum_{a,b} GW_{0,0}\big(F,H_{a}\big)\, GW_{0,0}\big(F,H_{b}\big)\, GW^{%
{\cal H}}_{S+dF,2}\big(H^{a},H^{b}\big)  \nonumber \\
\ +&\, \sum_{a}3\, \Big( GW^{{\cal H}}_{S+dF,0}\big(F,F,H_{a}\big)\, GW_{0,2}%
\big(\tau(H^{a})\big) \,+\, GW_{0,0}\big(F,F,H_{a}\big)\, GW^{{\cal H}%
}_{S+dF,2}\big(\tau(H^{a})\big)\,\Big)  \nonumber \\
\ -&\, \sum_{a,b}3\, GW^{{\cal H}}_{S+dF,0}\big(F,F,H_{a}\big)\, GW_{0,0}%
\big(H^{a},H_{b}\big)\, GW_{0,2}\big(H^{b}\big)  \nonumber \\
\ -&\, \sum_{a,b}3\, GW_{0,0}\big(F,F,H_{a}\big)\, GW^{{\cal H}}_{S+dF,0}%
\big(H^{a},H_{b}\big)\, GW_{0,2}\big(H^{b}\big)  \nonumber \\
\ -&\, \sum_{a,b}3\, GW_{0,0}\big(F,F,H_{a}\big)\, GW_{0,0}\big(H^{a},H_{b}%
\big)\, GW^{{\cal H}}_{S+dF,2}\big(H^{b}\big)  \nonumber \\
\ +&\, \sum_{a,b}\frac{13}{10}\, GW^{{\cal H}}_{S+dF,0}\big(F,F,H_{a},H_{b}%
\big)\, GW_{0,1}\big(H^{a}\big)\, GW_{0,1}\big(H^{b}\big)  \nonumber \\
\ +&\, \sum_{a,b}\frac{13}{10}\, GW_{0,0}\big(F,F,H_{a},H_{b}\big)\, GW^{%
{\cal H}}_{S+dF,1}\big(H^{a}\big)\, GW_{0,1}\big(H^{b}\big)  \nonumber \\
\ +&\, \sum_{a,b}\frac{13}{10}\, GW_{0,0}\big(F,F,H_{a},H_{b}\big)\, GW_{0,1}%
\big(H^{a}\big)\, GW^{{\cal H}}_{S+dF,1}\big(H^{b}\big)
\end{align*}

\begin{align}
\ +&\, \sum_{a,b}\frac{8}{5}\, GW^{{\cal H}}_{S+dF,1}\big(F,H_{a}\big)\,
GW_{0,0}\big(H^{a},F,H_{b}\big)\, GW_{0,1}\big(H^{b}\big)  \nonumber \\
\ +&\, \sum_{a,b}\frac{8}{5}\, GW_{0,1}\big(F,H_{a}\big)\, GW^{{\cal H}%
}_{S+dF,0}\big(H^{a},F,H_{b}\big)\, GW_{0,1}\big(H^{b}\big)  \nonumber \\
\ +&\, \sum_{a,b}\frac{8}{5}\, GW_{0,1}\big(F,H_{a}\big)\, GW_{0,0}\big(%
H^{a},F,H_{b}\big)\, GW^{{\cal H}}_{S+dF,1}\big(H^{b}\big)  \nonumber \\
\ -&\, \sum_{a,b}\frac{4}{5}\, GW^{{\cal H}}_{S+dF,0}\big(F,F,H_{a}\big)\,
GW_{0,1}\big(H^{a},H_{b}\big)\, GW_{0,1}\big(H^{b}\big)  \nonumber \\
\ -&\, \sum_{a,b}\frac{4}{5}\, GW_{0,0}\big(F,F,H_{a}\big)\, GW^{{\cal H}%
}_{S+dF,1}\big(H^{a},H_{b}\big)\, GW_{0,1}\big(H^{b}\big)  \nonumber \\
\ -&\, \sum_{a,b}\frac{4}{5}\, GW_{0,0}\big(F,F,H_{a}\big)\, GW_{0,1}\big(%
H^{a},H_{b}\big)\, GW^{{\cal H}}_{S+dF,1}\big(H^{b}\big)  \nonumber \\
\ +&\, \sum_{a,b}\frac{23}{240}\,\Big( GW^{{\cal H}}_{S+dF,0}\big(%
F,F,H_{a},H^{a},H_{b}\big)\, GW_{0,1}\big(H_{b}\big)\,+\, GW_{0,0}\big(%
F,F,H_{a},H^{a},H_{b}\big)\, GW^{{\cal H}}_{S+dF,1}\big(H_{b}\big)\,\Big)
\nonumber \\
\ +&\, \sum_{a,b}\frac{2}{48}\,\Big( GW^{{\cal H}}_{S+dF,0}\big(%
F,H_{a},H^{a},H_{b}\big)\, GW_{0,1}\big(H_{b},F\big)\,+\, GW_{0,0}\big(%
F,H_{a},H^{a},H_{b}\big)\, GW^{{\cal H}}_{S+dF,1}\big(H_{b},F\big)\,\Big)
\nonumber \\
\ -&\, \sum_{a,b}\frac{1}{80}\,\Big( GW^{{\cal H}}_{S+dF,1}\big(F,F,H_{a}%
\big)\, GW_{0,0}\big(H^{a},H_{b},H^{b}\big)\,+\, GW_{0,1}\big(F,F,H_{a}\big)%
\, GW^{{\cal H}}_{S+dF,0}\big(H^{a},H_{b},H^{b}\big)\,\Big)  \nonumber \\
\ +&\, \sum_{a,b}\frac{7}{30}\,\Big( GW^{{\cal H}}_{S+dF,0}\big(%
F,F,H_{a},H_{b}\big)\, GW_{0,1}\big(H^{a},H^{b}\big)\,+\, GW_{0,0}\big(%
F,F,H_{a},H_{b}\big)\, GW^{{\cal H}}_{S+dF,1}\big(H^{a},H^{b}\big)\,\Big)
\nonumber \\
\ +&\, \sum_{a,b}\frac{2}{30}\,\Big( GW^{{\cal H}}_{S+dF,0}\big(F,H_{a},H_{b}%
\big)\, GW_{0,1}\big(H^{a},H^{b},F\big)\,+\, GW_{0,0}\big(F,H_{a},H_{b}\big)%
\, GW^{{\cal H}}_{S+dF,1}\big(H^{a},H^{b},F\big)\,\Big)  \nonumber \\
\ -&\, \sum_{a,b}\frac{1}{30}\,\Big( GW^{{\cal H}}_{S+dF,0}\big(F,F,H_{a}%
\big)\, GW_{0,1}\big(H^{a},H_{b},H^{b}\big)\,+\, GW_{0,0}\big(F,F,H_{a}\big)%
\, GW^{{\cal H}}_{S+dF,1}\big(H^{a},H_{b},H^{b}\big)\,\Big)  \nonumber \\
\ +&\,\sum_{a,b}\frac{1}{576}\, GW^{{\cal H}}_{S+dF,0}\big(%
F,F,H_{a},H^{a},H_{b},H^{b}\big)  \label{TRR-g=2}
\end{align}
(cf. (17) of \cite{li}). Using the vanishing results in Lemma~\ref{GW=0},
one can simplify the right hand side of (\ref{TRR-g=2}) to have
\begin{align}  \label{pf-TRR}
GW_{S+dF,2}^{{\cal H}}\big(\tau_{1}(F),\tau_{2}(F)\big)\, &=\, \sum_{a,b} %
\big(-\frac{1}{80}\big)\, GW^{{\cal H}}_{S+dF,1}\big(F,F,H_{a}\big)\,
GW_{0,0}\big(H^{a},H_{b},H^{b}\big)  \nonumber \\
&+\,\sum_{a,b}\,\frac{1}{15}\, GW^{{\cal H}}_{S+dF,1}\big(F,H_{a},H_{b}\big)%
\, GW_{0,0}\big(F,H^{a},H^{b}\big)  \nonumber \\
&+\, \sum_{a,b}\,\frac{1}{576}\, GW^{{\cal H}}_{S+dF,0}\big(%
F,F,H_{a},H^{a},H_{b},H^{b}\big).
\end{align}
This can be further simplified by using Lemma~\ref{GW=0}\thinspace a. The
right-hand side of (\ref{pf-TRR}) becomes
\begin{equation}
-\frac{2}{3}\,GW_{S+dF,1}^{{\cal H}}\big(pt\big)\,+\,\sum_{a,b}\,\frac{1}{576%
}\,GW_{S+dF,0}^{{\cal H}}\big(F,F,H_{a},H^{a},H_{b},H^{b}\big).
\label{pf-TRR2}
\end{equation}
On the other hand, genus $g=0$ invariants with point constraints vanish by
dimensional reasons. This observation, combined with $\sum_{a}\,(H_{a}\cdot
A)(A\cdot H^{a})=A^{2}$, shows that
\begin{equation}
\sum_{a,b}\,\frac{1}{576}\,GW_{S+dF,0}^{{\cal H}}\big(%
F,F,H_{a},H^{a},H_{b},H^{b}\big)\,=\,\frac{4}{576}\,(4d-8)^{2}\,GW_{S+dF,0}^{%
{\cal H}}\,.  \label{pf-TRR3}
\end{equation}
Then, the proposition follows directly from (\ref{pf-TRR}), (\ref{pf-TRR2})
and (\ref{pf-TRR3}). \nopagebreak\hskip.1in {\ $\Box $ }\penalty10000 \hskip%
\parfillskip
\vskip1cm


\section{Proof of Theorem~\ref{T:Main}}

\bigskip

Our goal is to compute the genus $g=1$ family GW-invariants of K3 surfaces
for classes $A$ of index 2. By (\ref{T:BL}), it suffices to compute the
family GW-invariants of $E(2)$ for the classes $2(s+df)$. We introduce four
generating functions by the following formulas
\begin{eqnarray*}
M_{g}(t)\, &=&\,\sum \,GW_{2s+df,g}^{{\cal H}}\big(pt^{g}\big)\,t^{d}, \\
P_{g}(t)\, &=&\,\sum \,GW_{(s-3f)+d(2f),g}^{{\cal H}}\big(pt^{g}\big)\,t^{d},
\\
N_{g}(t)\, &=&\,\sum \,GW_{s+df,g}^{{\cal H}}\big(pt^{g}\big)\,t^{d}, \\
H_{g}(\,\cdot \,)(t)\, &=&\,\sum \,GW_{S+dF,g}^{{\cal H}}\big(\,\cdot \,\big)%
\,t^{d}.
\end{eqnarray*}
Notice that the coefficients of the even terms of $M_{1}\left( t\right) $
give the invariants $GW_{A,1}^{{\cal H}}$ for all index two classes $A$.
Therefore our main theorem \ref{T:Main} is equivalent to the following
proposition by restricting only to even terms.

\smallskip

\begin{prop}
The above generating functions satify the following relation,\
\[
M_{1}(t)\,=\,P_{1}(t)\,+\,2\,N_{1}(t^{2}).
\]
\end{prop}

\noindent {\bf Proof. } Since $(d-2)^{2}\,=\,d\,(d-1)\,-\,3\,d\,+\,4$, it
follows from Proposition~\ref{P:prop} that
\begin{align}
\hspace{-2cm}& \hspace{-0.5cm}GW_{S+dF,2}^{{\cal H}}\big(\tau _{1}(F),\tau
_{2}(F)\big)  \nonumber  \label{addition-1} \\
& \hspace{-0.5cm}=\,-\frac{2}{3}\,GW_{S+dF,1}^{{\cal H}}\big(pt\big)\,+\,%
\frac{1}{9}\,d\,(d-1)\,GW_{S+dF,0}^{{\cal H}}\,-\,\frac{1}{3}%
\,d\,GW_{S+dF,0}^{{\cal H}}\,+\,\frac{4}{9}\,GW_{S+dF,0}^{{\cal H}}.
\end{align}
Then, by combining (\ref{addition-1}) and the definition of $H_{g}(\,\cdot
\,)(t)$, we obtain
\begin{equation}
H_{2}\big(\tau _{1}(F),\tau _{2}(F)\big)(t)\,=\,\frac{1}{9}\,H_{0}^{\prime
\prime }(t)\,-\,\frac{1}{3}\,t\,H_{0}^{\prime }(t)\,+\,\frac{4}{9}%
\,H_{0}(t)-\,\frac{2}{3}\,H_{1}(t).  \label{TRR}
\end{equation}
In \cite{ll}, we used the symplectic sum formula of \cite{ip} to obtain
\begin{align}
& H_{2}\big(\tau _{1}(F),\tau _{2}(F)\big)(t)\,-2\,H_{1}(t)  \nonumber
\label{ODE-1} \\
& =\,\frac{20}{3}\,G_{2}(t)\,t\,H_{0}^{\prime }(t)\,-\,\left(
64\,G_{2}^{2}(t)\,+\,\frac{40}{3}\,G_{2}(t)\,-\,8\,t\,G_{2}^{\prime
}(t)\right) H_{0}(t).
\end{align}
The equations (\ref{ODE-1}) and (\ref{TRR}) then yield
\begin{align}
& \frac{1}{9}\,t^{2}\,H_{0}^{\prime \prime }(t)\,-\,\frac{1}{3}%
\,t\,H_{0}^{\prime }(t)\,+\,\frac{4}{9}\,H_{0}(t)\,-\,\frac{8}{3}\,H_{1}(t)
\nonumber  \label{ODE-2} \\
& =\,\frac{20}{3}\,G_{2}(t)\,t\,H_{0}^{\prime }(t)\,-\,\left(
64\,G_{2}^{2}(t)\,+\,\frac{40}{3}\,G_{2}(t)\,-\,8\,t\,G_{2}^{\prime
}(t)\right) H_{0}(t).
\end{align}
On the other hand, we have
\begin{align}
\frac{1}{8}\,N_{0}(t^{2})\,& =\,M_{0}(t)-P_{0}(t),  \label{main-ll} \\
t\,\frac{d}{dt}\,N_{0}(t)\,& =\,24\,G_{2}(t)\,N_{0}(t)\,+\,N_{0}(t),
\label{pri-g=0} \\
N_{1}(t)\,& =\,\Big(t\,\frac{d}{dt}\,G_{2}(t)\Big)N_{0}(t),  \label{pri-g=1}
\end{align}
(cf. \cite{ll,bl1,l2}). It then follows from (\ref{main-ll}) and (\ref
{pri-g=0}) that
\begin{equation}
\frac{d}{dt}\,\big(\,M_{0}(t)\,-\,P_{0}(t)\,\big)\,=\,\big(%
\,48\,G_{2}(t^{2})\,+\,2\,\big)\,\big(\,M_{0}(t)\,-\,P_{0}(t)\,\big).
\label{bas-res}
\end{equation}
Recalling $(S,F)$ is either $(2s,f)$ or $(s-3f,2f)$, we combine (\ref{ODE-2}%
) and (\ref{bas-res}) to obtain
\begin{equation}
M_{1}(t)\,-\,P_{1}(t)\,=\,\big(\,4\,t^{2}\,G_{2}^{\prime
}(t^{2})\,+\,3\,F(t)\,\big)\big(\,M_{0}(t)\,-\,P_{0}(t)\,\big)
\label{com-com}
\end{equation}
where $F(t)\,=\,32\,G_{2}^{2}(t^{2})\,-\,40\,G_{2}(t^{2})\,G_{2}(t)\,+\,8%
\,G_{2}^{2}(t)\,-\,t\,G_{2}^{\prime }(t).$ Then, we have
\[
M_{1}(t)\,-\,P_{1}(t)\,=\,16\,t^{2}\,G_{2}^{\prime }(t^{2})\,\big(%
\,M_{0}(t)\,-\,P_{0}(t)\,)\,=\,2\,t^{2}\,G_{2}^{\prime
}(t^{2})\,N_{0}(t^{2})\,=\,2\,N_{1}(t^{2})
\]
where the first equality follows from (\ref{com-com}) and the fact $%
F(t)\,=\,4\,t^{2}\,G_{2}^{\prime }(t^{2})$ \cite{ll}, the second equality
follows from (\ref{main-ll}), and the last equality follows from (\ref
{pri-g=1}). \nopagebreak\hskip.1in {\ $\Box $ }\penalty10000 \hskip%
\parfillskip

\vskip1cm


\bigskip


\begin{thebibliography}{BL1}



\bibitem[BL1]{bl1}  J. Bryan and N.C. Leung, {\em The enumerative geometry
of K3 surfaces and modular forms}, J. Amer. Math. Soc. {\bf 13} (2000),
371-410.












\bibitem[BL2]{bl2}  J. Bryan and N.C. Leung, {\em Counting curves on
irrational surfaces}, Survey of Differential Geometry. {\bf 5} (1999),
313-339.

\bibitem[Ge]{g}  E. Getzler, {\em Topological recursion relations in genus 2}%
, In ''Integrable systems and algebraic geometry (Kobe/Kyoto, 1997).'' World
Sci. Publishing, River Edge, NJ, 198, pp 73-106.

\bibitem[G]{go}  L. G\"{o}ttsche, {\em A conjectural generating function for
numbers of curves on surfaces}, preprint, alg-geom/9711012






\bibitem[IP]{ip}  E. Ionel and T. Parker, {\em The Symplectic Sum Formula
for Gromov-Witten Invariants}, to appear in Ann. Math



\bibitem[KM]{km1}  M. Kontsevich and Y.I. Manin, {\em Relations between the
correlators of the topological sigma model coupled to gravity}, Commun.
Math. Phys. {\bf 196} (1998), 385-398.




\bibitem[L1]{l1}  J. Lee, {\em Family Gromov-Witten Invariants for
K\"{a}hler Surfaces}, to appear in Duke Math. J.

\bibitem[L2]{l2}  J. Lee, {\em Counting Curves in Elliptic Surfaces by
Symplectic Methods}, preprint, math.SG/0307358.


\bibitem[Li]{li}  X. Liu, {\em Quantum product on the big phase space and
the Virasoro conjecture}, Adv. math. {\bf 169} (2002), 313-375.

\bibitem[LL]{ll}  J. Lee and N.C. Leung, {\em Yau-Zaslow formula on K3
surfaces for non-primitive classes}, preprint, math.SG/0404537.




\bibitem[LT]{lt}  J. Li and G. Tian, {\em Virtual moduli cycles and
Gromov-Witten invariants of general symplectic manifolds}, Topics in
symplectic $4$-manifolds (Irvine, CA, 1996), 47--83, First Int. Press Lect.
Ser., I, International Press, Cambridge, MA, 1998.


\bibitem[M]{m}  D. Mumford, {\em Towards an enumerating geometry of the
moduli space of curves}, in Arithmetic and Geometry, M. Artin and J. Tate,
eds., Birkh\"{a}user, 1995, 401-417.







\bibitem[RT]{rt}  Y. Ruan and G. Tian, {\em Higher genus symplectic
invariants and sigma models coupled with gravity}, Invent. Math. {\bf 130}
(1997), 455-516.





\bibitem[YZ]{yz}  S.T. Yau and E. Zaslow, {\em BPS States, String Duality,
and Nodal Curves on K3}, Nuclear Phys. B {\bf 471} (1996), 503-512.
\end{thebibliography}
\end{document}